\documentclass[12pt,a4paper,draft]{article}
\usepackage[T1]{fontenc}
\usepackage{amsmath}
\usepackage{amssymb}
\usepackage{mathrsfs}
\newtheorem{thm}{Theorem}[section]
\newtheorem{lem}[thm]{Lemma}
\newtheorem{cor}[thm]{Corollary}
\def\fd{\hfill$\square$\par\noindent\vskip0.2cm}
\newcommand{\pre}{\textbf{Proof}}
\newcommand{\preu}{\textbf{Proof of Theorem }}
\newcommand{\R}{\mathbb{R}}
\newcommand{\C}{\mathbb{C}}

\newcommand{\N}{\mathbb{N}}
\newcommand{\Z}{\mathbb{Z}}
\newcommand{\supp}{\mbox{\rm supp}}

\newcommand{\im}{\mbox{\rm Im}}
\numberwithin{equation}{section}
\begin{document}
\title{On covariant functions and distributions under the action of a compact group}
\author {Anouar SAIDI \footnote{\scriptsize  Faculty of Sciences of
 Monastir. Department
of Mathematics. 5019 Monastir.Tunisia.\vfill 
Email: saidi.anouar@yahoo.fr 
}}
\maketitle
\begin{scriptsize}
\begin{center}
\begin{minipage}{11cm}
{\bf Abstract. } 
Let $G$ be a compact subgroup of $GL_n(\R)$ acting linearly on a finite dimensional vector space $E$. B. Malgrange has shown that the space $\mathcal{C}^\infty(\R^n,E)^G$ of $\mathcal{C}^\infty$ and $G$-covariant functions is a finite module over the ring $\mathcal{C}^\infty(\R^n)^G$ of $\mathcal{C}^\infty$ and $G$-invariant functions. First, we generalize this result for the Schwartz space $\mathscr{S}(\R^n,E)^G$ of $G$-covariant functions. Secondly, we prove that any $G$-covariant distribution can be decomposed into a sum of $G$-invariant distributions multiplied with a fixed family of $G$-covariant polynomials. This gives a generalization of  an Oksak result proved in ([O]). 
\end{minipage}
\end{center}


\end{scriptsize}
\section{Introduction}
Let $G$ be a compact Lie group acting linearly on $\R^n$ and $\rho$ a representation  of $G$ in a finite dimensional complex vector space $E$. 

Let $\mathscr{A}(\R^n,E)$ be either of the spaces $\mathcal{C}^\infty(\R^n,E)$, $\mathcal{C}_c^\infty(\R^n,E)$ or $\mathscr{S}(\R^n,E)$. 
A function $f$ in $\mathscr{A}(\R^n,E)$ is said to be $G$-covariant if it verifies 
$$\rho(g)^{-1}\circ f\circ g=f,\quad\textrm{ for all }g\in G.$$
We denote by $\mathscr{A}(\R^n,E)^G$ the space of $G$-covariant functions and by $\mathscr{A}(\R^n)^G$ the space of $G$-invariant functions. 

Let $E^*$ be the dual space of $E$ and $\rho^*$ the adjoint representation of $\rho$ in $E^*$. We consider $\mathscr{A}'(\R^n,E)$ the dual space of $\mathscr{A}(\R^n,E^\ast)$.  
A distribution $T$ in $\mathscr{A}'(\R^n,E)$ is said to be $G$-covariant if it verifies  $$\langle T,\rho(g)^{-1}\circ f\circ g\rangle=\langle T,f\rangle,\quad\textrm{ for all }f \in\mathscr{A}(\R^n,E^*) \textrm{ and }g\in G.$$ 
Let $\mathscr{A}'(\R^n,E)^G$ be the space  of $G$-covariant distributions and  $\mathscr{A}'(\R^n)^G$ the space of $G$-invariant distributions.

Representations of $G$-covariant functions and distributions in terms of $G$-covariant polynomials are of frequent use in theoretical physics. B. Malgrange (see \cite{POE} of V. Poénaru) has obtained such  representation for functions in $\mathcal{C}^\infty(\R^n,E)^G$.
More precisely, he has proved that there exists a family of $G$-covariant polynomials $P_1,\ldots,P_k$ such that every function $f$ in $\mathcal{C}^\infty(\R^n,E)^G$ can be represented in the form \begin{equation}\label{eq01}f=f_1P_1+\ldots+f_rP_k,\quad\textrm{ where }  f_1,\ldots,f_k\in\mathcal{C}^\infty(\R^n)^G.\end{equation} 
 
A.I. Oksak in \cite{OK} has considered the case when the subspace $\mathscr{P}(\R^n,E)^G$ of $G$-covariant polynomials  is a free module over the ring $\mathscr{P}(\R^n)^G$ of $G$-invariant polynomials. This is equivalent to say that the map $I:(\mathscr{P}(\R^n)^G)^k\longrightarrow\mathscr{P}(\R^n,E)^G$ defined by 
$$(p_1,\ldots,p_k)\longmapsto\sum_{i=1}^kp_i\;P_i,$$ is an isomorphism.
Under this assumption, he established that the map $\Phi:(\mathscr{A}(\R^n)^G)^k\longrightarrow\mathscr{A}(\R^n,E)^G$ defined by 
$$(f_1,\ldots,f_k)\longmapsto\sum_{i=1}^kf_i\;P_i,$$ is a topological  isomorphism when the vector space 
 $\mathscr{A}(\R^n,E)^G$ is equal to $\mathcal{C}^\infty(\R^n,E)^G$, $\mathcal{C}_c^\infty(\R^n,E)^G$ or $\mathscr{S}(\R^n,E)^G$.

A.I. Oksak has also studied the representation of $G$-covariant distributions. He has proved, under the above assumption, that every distribution $T$ in $\mathscr{A}'(\R^n,E)^G$ can be represented in the form \begin{equation}\label{eq02}
T=\sum_{i=1}^kT_i\;P_i,
\quad\textrm{ where } T_1,\ldots,T_k\in\mathscr{A}'(\R^n)^G.\end{equation}

Our aim in this  paper is to generalize the results of B. Malgrange  and A.I. Oksak. In section 2, we prove the analogue of the result of B. Malgrange in the spaces $\mathcal{C}^\infty_c(\R^n,E)^G$ and $\mathscr{S}(\R^n,E)^G$. In section 3, we  establish the decomposition (\ref{eq02})   without any assumption on the action of the compact group $G$.

\section{Representation of $G$-covariant functions}
Let $G$ be a compact Lie group and $\mu$ its normalized Haar measure. We suppose that $G$ acts linearly on $\R^n$ and we fix  an  euclidean $G$-invariant norm on $\R^n$. Let $\rho$ be a representation of $G$ on a finite dimensional complex vector space $E$. We choose a $G$-invariant norm  $||\;\;||_E$ on $E$. By duality, $G$ acts on $E^\ast$ the dual space of $E$. We consider $\mathcal{C}_c^\infty(\R^n,E)$ the space of compact support functions in $\mathcal{C}^\infty(\R^n,E)$. For a compact $K$ in $\mathbb{R}^n$,   $\mathcal{C}_K^\infty(\R^n,E)$ denotes
the space of  functions in $\mathcal{C}_c^\infty(\R^n,E)$  such that $\supp(f)\subset K$  and $\mathscr{S}(\R^n,E)$ the space of functions $f$ in $\mathcal{C}^\infty(\R^n,E)$ satisfying 
$$\sup_{x\in\;\R^n}(1+||x||^2)^m\;||D^\alpha f(x)||_E<\infty,$$ for all $m\in\N$ and $\alpha\in\N^n$. 

We equip   the spaces $\mathcal{C}^\infty(\R^n,E)$, $\mathcal{C}_c^\infty(\R^n,E)$ and $\mathscr{S}(\R^n,E)$  with their usual topologies. For $g\in G$ and $f\in\mathscr{A}(\R^n,E)$, we denote by $g\hskip-0,01cm\cdot\hskip-0,01cm f$ the function $\rho(g)^{-1}\circ f\circ g$ which also belongs to $\mathscr{A}(\R^n,E)$. We say that  a function $f$ in $\mathscr{A}(\R^n,E)$ is $G$-covariant if it verifies $$g\hskip-0,02cm\cdot\hskip-0,02cm f=f,\quad\textrm{ for all }g\in G.$$ 

It is clear that $\mathscr{A}(\R^n,E)^G$ is an $\mathscr{A}(\R^n)^G$-module. Our aim in this section is to prove that this module is generated by a finite family of $G$-covariant polynomials. 
The case $\mathcal{C}^\infty(\R^n,E)$, is already established by B. Malgrange (see \cite{POE}) and the same proof is valid to establish that the set of $G$-covariant polynomials is a module of finite type over the ring of $G$-invariant polynomials.
We fix $P_1,\ldots,P_k$ a family of $G$-covariant homogeneous polynomials which generate this module.
So we have :
\begin{lem}\label{lem0} 
For all $f\in\mathcal{C}^\infty(\R^n,E)^G$ there exists a  family of functions $f_1,\ldots,f_k$ in $\mathcal{C}^\infty(\R^n)^G$ such that $$f=\sum_{i=1}^k f_i\;P_i.$$
\end{lem}
As corollary of the above lemma, we have 
\begin{cor}\label{cor1} 
For all $f\in\mathcal{C}_c^\infty(\R^n,E)^G$ there exists a  family of functions $f_1,\ldots,f_k$ in $\mathcal{C}_c^\infty(\R^n)^G$ such that $$f=\sum_{i=1}^k f_i\;P_i.$$
\end{cor}
The main result of this section is :
\begin{thm}\label{th2}
For all $f\in\mathscr{S}(\R^n,E)^G$ there exists a  family of functions $f_1,\ldots,f_k$ in $\mathscr{S}(\R^n)^G$ such that $$f=\sum_{i=1}^k f_i\;P_i.$$
\end{thm} 
To prove this theorem, we need the following lemmas.

\vskip0,3cm
For all $0\leq r<R$, we consider the sets $B(0,r)=\{x\in\R^n,\;||x||<r\}$, $S(0,r)=\{x\in\R^n,\;||x||=r\}$, $\overline{B}(0,r)=B(0,r)\cup S(0,r)$ and $C(r,R)=\{x\in\R^n,\;r\leq||x||\leq R\}.$ 

As an improvement of Corollary \ref{cor1}, we have
\begin{lem}\label{lem3}
Let $f$ be in $\mathcal{C}_{\overline{B}(0,1)}^\infty(\R^n,E)^G$. There exists a family of functions $f_1,\ldots,f_k$ in $\mathcal{C}_{\overline{B}(0,1)}^\infty(\R^n)^G$ such that    $$f=\sum_{i=1}^k f_i\;P_i.$$ 
\end{lem}
\pre\\
Let $f$ be in $\mathcal{C}_{\overline{B}(0,1)}^\infty(\R^n,E)^G.$ By Corollary \ref{cor1} we have $$f=\sum_{i=1}^k{f_i\;P_i},\quad\textrm{ with } f_1,\ldots,f_k\in \mathcal{C}^\infty_c(\R^n)^G.$$ We have to prove that the functions $f_1,\ldots,f_k$ can be chosen with supports in $\overline{B}(0,1)$. 
For this, we need to use Lemma 8.1.1 of \cite{BO1}. This lemma treats the case $E=\C$ but it is also true for any finite dimensional vector space. 

Let $\gamma$ be in $\mathcal{C}^\infty(\R^n)^G $ satisfying $\gamma\equiv1$ on $\overline{B}(0,\frac{1}{2})$ and $\supp(\gamma)\subset\overline{B}(0,\frac{2}{3}).$ Then $f=\gamma f+(1-\gamma)f.$ 

The function $(1-\gamma)\;f\in\mathcal{C}_{C(\frac{1}{2},1)}^\infty(\R^n,E)^G
$, so Lemma 8.1.1 of \cite{BO1} implies that  $$(1-\gamma(x))\;f(x)=\sum_{j=1}^r g_j(||x||^2)\;
 \varphi_{j}(x),\quad\textrm{ for all } x\in\R^n,$$
  with $g_j\in\mathcal{C}^\infty_{[\frac{1}{4},1]}(\R)$ and $\varphi_{j}\in\mathcal{C}_{C(\frac{1}{2},1)}^\infty(\R^n,E).$ Since $G$ is a compact group, we can suppose that the functions $\varphi_j$ are $G$-covariant. Then by Corollary \ref{cor1}, we have  $$\varphi_{j}=\sum_{i=1}^k\varphi_{ij}\;P_i,\quad\textrm{ for all }1\leq j\leq r,$$ where $\varphi_{ij}\in\mathcal{C}_c^\infty(\R^n)^G$. So we get    $$(1-\gamma(x))\;f(x)=\sum_{i=1}^k\bigg(\sum_{j=1}^rg_j(||x||^2)\varphi_{ij}(x)\bigg)P_i(x),\quad\textrm{ for all } x\in\R^n.$$ The result follows for the function $(1-\gamma)f$ since the functions $g_j(||\;.\;||^2)\;\varphi_{ij}$ are in $\mathcal{C}_{\overline{B}(0,1)}^\infty(\R^n)^G$. 
 
According to Corollary \ref{cor1}, there exist $h_1,\ldots,h_k\in \mathcal{C}^\infty_c(\R^n)^G$ such that $$\gamma f=\sum_{i=1}^k{h_i\;P_i}.$$
Let $\theta$ be in $\mathcal{C}_c^\infty(\R^n)^G$ satisfying $\theta\equiv1$ on $\overline{B}(0,\frac{2}{3})$ and $\supp(\theta)\subset\overline{B}(0,1).$ It is clear that 
$$\gamma f=\theta\;\gamma\;f=\sum_{i=1}^k{(\theta\;h_i)\;P_i},$$ with  $\theta\; h_i\in \mathcal{C}_{\overline{B}(0,1)}^\infty(\R^n)^G$. This completes the proof.\fd
\begin{lem}\label{lem1}
Let $h$ be in $\mathcal{C}^\infty(\R^n,E)^G$ such that for all $\alpha$ in $\N^n$, $D^\alpha h\equiv0$ on $S(0,1)$. Then, for all $p$ in $\N$,  \begin{displaymath}
\lim_{||y||\to1^-}\frac{h(y)}{(1-||y||^2)^p}=0.\end{displaymath}
\end{lem}
\pre\\  
Let $h$ be in $\mathcal{C}^\infty(\R^n,E)^G$ and $p$ in $\N$. We set $M=\displaystyle\sup_{z\in\overline{B}(0,2)}||D^{p+1}h(z)||$. By Taylor formula we have 
$$h(y)=\int_0^1{\frac{(1-t)^p}{p!}\; D^{p+1}h(y_0+t(y-y_0))\cdot(y-y_0)^{(p+1)}}\;dt,$$ for all $y\in\overline{B}(0,1)$ and $y_0$ in $S(0,1)$. Then, $$|h(y)|\leq M\;||y-y_0||^{p+1}, \quad\textrm{ for all } y\in\overline{B}(0,1) \textrm{ and } y_0\in S(0,1).$$
Let $y\in B(0,1)\setminus\{0\}$ and $y_0=\frac{y}{\left\|y\right\|}$. Then  $\left\|y-y_0\right\|=1-\left\|y\right\|$ and 
$$ \frac{|h(y)|}{(1-||y||^2)^p}=
\frac{|h(y)|}{(1+||y||)^p||y-y_0||^p}\leq M\;||y-y_0||.$$
This gives the result since $$\lim_{||y||\to1^-}||y-y_0||=0.$$\fd

The function  $\varphi:x\longmapsto\frac{x}{\sqrt{1-\left\|x\right\|^2}}$  
is a bijection from $B(0,1)$ to $\R^n$ and its inverse  function  is $$\psi:x\longmapsto\frac{x}{\sqrt{1+||x||^2}}.$$

For $f$ in $\mathscr{S}(\R^n,E)$, we denote by $I(f)$ the function on $\R^n$ defined by :\begin{displaymath} I(f)(x)=\left\{ \begin{array}{ll}
f\circ\varphi(x)&\textrm{if $x\in B(0,1)$}\\
    0 &  \textrm{ otherwise }
    \end{array} \right..
\end{displaymath} 
The map $I$ is one to one and satisfies :
\begin{lem}\label{lem2}
For all $f\in\mathscr{S}(\R^n,E)$, the function $I(f)$ is in $\mathcal{C}_{\overline{B}(0,1)}^\infty(\R^n,E)$ and the map $I:\mathscr{S}(\R^n,E)\longrightarrow\mathcal{C}_{\overline{B}(0,1)}^\infty(\R^n,E)$ is a topological isomorphism.
 \end{lem}
 \pre\\
Let $f$ be in $\mathscr{S}(\R^n,E)$. To prove that $I(f)$ is $\mathcal{C}^\infty$ it is enough to show that  $$ \lim_{\left\|x\right\|\to 1^-}D^\alpha(f\circ\varphi)(x)=0,\quad\textrm{ for all } \alpha\in\N^n.$$
It is clear that, for all $x\in B(0,1)$, $D^\alpha(f\circ\varphi)(x)$ can be written as a sum of elements of the form 
$$D^\gamma f(\varphi(x))\;D^{\lambda_1}\varphi_{i_1}(x)\cdots D^{\lambda_r}
\varphi_{i_r}(x),$$ with $\lambda_i, \gamma\in\N^n$,  and $$ D^{\lambda_j}\varphi_{i_j}(x)=\frac{Q_{i,j}(x)}{(1-||x||^2)^s},$$  where  $Q_{i,j}$ is a polynomial and  $s\in\frac{1}{2}\;\!\N$. Since  $$1-||x||^2=\frac{1}{1+||\varphi(x)||^2},$$  
then there exist a polynomial $P$ and  $s'\in\frac{1}{2}\;\!\N$ such that 
\begin{equation}\label{eq1} D^\gamma f(\varphi(x))\;D^{\lambda_1}\varphi_{i_1}(x)\cdots D^{\lambda_r}\varphi_{i_r}(x)=P(x)\;(1+||\varphi(x)||^2)^{s'}\;D^\gamma f(\varphi(x)).\end{equation} 
As $f\in\mathscr{S}(\R^n,E)$ and $\lim_{||x||\to 1^-}||\varphi(x)||=+\infty,$  then \begin{displaymath}\lim_{||x||\to 1^-}D^\alpha(f\circ\varphi)(x)=0,\end{displaymath}
and we deduce that $I(f)\in\mathcal{C}_{\overline{B}(0,1)}^\infty(\R^n,E)$. 

To establish that $I$ is onto we will prove that if $g\in \mathcal{C}_{\overline{B}(0,1)}^\infty(\R^n,E)$, then $g\circ\psi\in\mathscr{S}(\R^n,E)$ and $I(g\circ \psi)=g$. We fix $r\in\N$ and $\beta\in\N^n$. As in (\ref{eq1}) there exist $m\in\Z$ and $c>0$ such that 
$$(1+||x||^2)^r||D^\beta(g\circ\psi)(x)||_E\leq c\;(1+||x||^2)^m\sup_{\gamma\leq\beta}||D^\gamma g(\psi(x))||_E.$$ Since $\displaystyle\lim_{||x||\to+\infty}||\psi(x)||=1$ and $\displaystyle 1+||x||^2=\frac{1}{1-||\psi(x)||^2}$, then Lemma \ref{lem1} insures that 
$$\lim_{||x||\to+\infty}(1+||x||^2)^m||D^\gamma g(\psi(x))||_E=0,$$ and we deduce that $$\lim_{||x||\to+\infty}(1+||x||^2)^r||D^\beta(g\circ\psi)(x))||_E=0.$$ Hence $g\circ\psi\in\mathscr{S}(\R^n,E)$, and it is easy to verify that $I(g\circ\psi)=g$. 

It follows from the closed graph theorem that $I$ is a topological isomorphism.\fd
\noindent\preu{\bf\ref{th2}}\\
Let $f$ be in $\mathscr{S}(\R^n,E)^G$. Since we have 
$$\varphi(g\;x)=g.(\varphi\;(x)),\quad\textrm{ for all }g\in G \textrm{ and }x\in B(0,1),$$
then the function $I(f)$ is in $\mathcal{C}_{\overline{B}(0,1)}^\infty(\R^n,E)^G$ and by Lemma \ref{lem3}, we have  $$I(f)(x)=\sum_{i=1}^kg_i(x)\;P_i(x),\quad\textrm{ for all } x\in\R^n,$$ with  $g_1,\ldots,g_k\in\mathcal{C}_{\overline{B}(0,1)}^\infty(\R^n)^G.$
So for all $y$ in $\R^n$, we have  $$f(y)=I(f)(\psi(y))=\sum_{i=1}^kg_i(\psi(y))\;P_i(\psi(y)).$$  Let $s_1,\ldots,s_k$ be the degrees of homogeneity of $P_1,\ldots,P_k$. Then $$P_i(\psi(y))=P_i\big(\frac{y}{\sqrt{1+||y||^2}}\big)=\big(\frac{1}{\sqrt{1+||y||^2}}\big)^{s_i}P_i(y),\quad\textrm{ for all }1\leq i\leq k.$$ Hence 
$$f(y)=\sum_{i=1}^kf_i(y)\;P_i(y),$$
with $$f_i(y)=\big(\frac{1}{\sqrt{1+||y||^2}}\big)^{s_i}\;g_i\big(\frac{y}{\sqrt{1+||y||^2}}\big)\in\mathscr{S}(\R^n)^G.$$ 
To conclude, we have to prove that $f_1,\ldots,f_k\in\mathscr{S}(\R^n)$.
We fix $1\leq i\leq k$ and we consider the function $h_i:\R^n\longrightarrow\C$ defined by  
\begin{displaymath} h_i(x)=\left\{ \begin{array}{ll}
(\sqrt{1-||x||^2})^{s_i}\;g_i(x)&\textrm{ if } x\in B(0,1)\\
    0 &  \textrm{ otherwise }
    \end{array} \right..
\end{displaymath} 
The partial derivatives of $h_i$ in $B(0,1)$ are of the form $P(x)(1-||x||^2)^r D^\gamma g_i(x)$ with $r\in\frac{1}{2}\:\!\Z$, $\gamma\in\N^n$ and $P$ a polynomial. Hence, by Lemma \ref{lem1}, we conclude that    $\displaystyle\lim_{||x||\to 1^-}D^\alpha h_i(x)=0$ for all $\alpha\in\N^n$ and  so $h_i\in \mathcal{C}_{\overline{B}(0,1)}^\infty(\R^n)$. Finally, we  verify that $f_i=h_i\circ\psi=I^{-1}(h_i)$ and that $f_i$ is $G$-invariant. This completes the proof. \fd
\section{Representation of $G$-covariant distributions}
Let  $\mathscr{A}'(\R^n,E)$  be the topological dual of the space $\mathscr{A}(\R^n,E^*)$ equipped with its weak topology.  
A distribution $T$ in $\mathscr{A}'(\R^n,E)$ is said to be $G$-covariant if it verifies  \begin{equation}\label{eq2}\langle T,g\hskip-0,009cm\cdot\hskip-0,01cm  f\rangle=\langle T,f\rangle,\quad\textrm{ for all }f \in\mathscr{A}(\R^n,E^*) \textrm{ and }g\in G.\end{equation} We denote by $\mathscr{A}'(\R^n,E)^G$ the space of $G$-covariant distributions. In the particular case of the trivial representation, a $G$-covariant distribution $T$ is $G$-invariant and we denote by $\mathscr{A}'(\R^n)^G$ the space $\mathscr{A}'(\R^n,\C)^G$.

For $f\in\mathscr{A}(\R^n,E^*)$ and $h\in\mathscr{A}(\R^n,E)$, we define the function $\langle f,h\rangle$ in $\mathscr{A}(\R^n)$ by $$\langle f,h\rangle(x)=\langle f(x),h(x)\rangle,\quad\textrm{ for all } x\in\R^n.$$ 

Let $\theta\in\mathscr{A}'(\R^n)$ and $h\in\mathscr{A}(\R^n,E)$. We denote by $\theta\,h$ the distribution in $\mathscr{A}'(\R^n,E)$ defined by  $$\langle\theta\;h,f\rangle=\langle\theta,\langle f,h\rangle\rangle,\quad\textrm{ for all }f\in\mathscr{A}(\R^n,E^*).$$ If  $\theta\in\mathscr{A}'(\R^n)^G$ and $h\in\mathscr{A}(\R^n,E)^G$, then  $\theta\,h\in\mathscr{A}'(\R^n,E)^G$.

Then we have the following theorem :
\begin{thm}\label{th3}
For all $\theta$ in $\mathscr{A}'(\R^n,E)^G$ there exists a family of distributions $\theta_1,\ldots,\theta_k$ in $\mathscr{A}'(\R^n)^G$ such that  $\displaystyle\theta=\sum_{i=1}^k\theta_i\;P_i.$ 
        \end{thm}
Before the proof of this theorem, we establish  some lemmas.

\vskip 0.3cm
For $\varphi\in\mathcal{C}^{\infty}_c(\R^n,E)$, let $T_\varphi$ be the distribution in $\mathscr{A}'(\R^n,E)$ defined by $\displaystyle T_\varphi(f)\!=\int_{\R^n}\langle\varphi(x),f(x)\rangle\:dx$. Then we can identify $\mathcal{C}^{\infty}_c(\R^n,E)$ with a subspace of  $\mathscr{A}'(\R^n,E)$. Moreover, if $\varphi$ is a $G$-covariant function, then $T_\varphi$ is a $G$-covariant distribution and we have the following result.
\begin{lem}\label{lem4}
The space $\mathcal{C}^{\infty}_c(\R^n,E)^G$ is dense in $\mathscr{A}'(\R^n,E)^G$.
\end{lem}
\pre\\
Let $f$ be in $\mathscr{A}(\R^n,E^\ast)$. The map $g\longmapsto g\hskip-0,005cm\cdot\hskip-0,007cm  f$ from $G$ to $\mathscr{A}(\R^n,E^\ast)$ is clearly continuous. By Proposition 5 in \cite{BO2}, there exists a unique function $\mathscr{M}(f)$ in $\mathscr{A}(\R^n,E^\ast)$, such that for all $T\in\mathscr{A}'(\R^n,E)$, we have 
\begin{equation}\label{eq4}\langle T,\mathscr{M}(f)\rangle=\int_G\langle T, g\hskip-0,005cm\cdot\hskip-0,007cm  f\rangle\;d\mu(g).\end{equation}  It follows from this equality that \begin{equation}\label{eq5}\mathscr{M}(f)=\mathscr{M}(g\hskip-0,005cm\cdot\hskip-0,007cm  f),\quad\textrm{ for all } f\in\mathscr{A}(\R^n,E^\ast) \textrm{ and all }g\in G.\end{equation} 
For $T=\delta_x$, the Dirac distribution on $x\in\R^n$, the formula (\ref{eq4}) gives $$\mathscr{M}(f)(x)=\int_G g\hskip-0,005cm\cdot\hskip-0,007cm  f\;(x)\;d\mu(g).$$
Then we can verify that :

\vskip0,2cm
i)\; $\mathscr{M}(f)\in\mathscr{A}(\R^n,E^\ast)^G$. 

\vskip0,2cm 
ii) \;$\mathscr{M}(f)=f$ if $f$ is $G$-covariant.

\vskip0,2cm
iii)\; The linear map $\mathscr{M}:\mathscr{A}(\R^n,E^\ast)\longrightarrow\mathscr{A}(\R^n,E^\ast)$ is continuous (according to the graph theorem).

\vskip0,2cm
iv)\;  $\mathscr{M}(\mathcal{C}^{\infty}_c(\R^n,E))=\mathcal{C}^{\infty}_c(\R^n,E)^G$. 

\vskip0,3cm
Let  ${}^t\hskip-0,12cm\mathscr{M}$ be the transpose map of $\mathscr{M}$. It is continuous for the weak topologies and verifies 
\begin{equation}\label{eq6}{}^t\hskip-0,12cm\mathscr{M}(T)=T,\quad\textrm{ for all } T\in\mathscr{A}'(\R^n,E)^G,\end{equation} and  \begin{equation}\label{eq7}{}^t\hskip-0,12cm\mathscr{M}(T_\varphi)=T_{\mathscr{M}(\varphi)},\quad\textrm{ for all }\varphi\in\mathcal{C}^{\infty}_c(\R^n,E).\end{equation}

 Since $\mathcal{C}^{\infty}_c(\R^n,E)$ is dense in  $\mathscr{A}'(\R^n,E)$, then ${}^t\hskip-0,12cm\mathscr{M}(\mathcal{C}^{\infty}_c(\R^n,E))$   is dense in ${}^t\hskip-0,12cm\mathscr{M}(\mathscr{A}'(\R^n,E))$. 
 
 To establish the lemma we will prove that ${}^t\hskip-0,12cm\mathscr{M}(\mathscr{A}'(\R^n,E))=\mathscr{A}'(\R^n,E)^G$ and  ${}^t\hskip-0,12cm\mathscr{M}(\mathcal{C}^{\infty}_c(\R^n,E))=\mathcal{C}^{\infty}_c(\R^n,E)^G$.

 The inclusion ${}^t\hskip-0,12cm\mathscr{M}(\mathscr{A}'(\R^n,E))\subset\mathscr{A}'(\R^n,E)^G$ follows from the formula (\ref{eq5}) and the formula (\ref{eq6}) gives the inverse inclusion. 
Finally, the equality ${}^t\hskip-0,12cm\mathscr{M}(\mathcal{C}^{\infty}_c(\R^n,E))=\mathcal{C}^{\infty}_c(\R^n,E)^G$ follows from the formula (\ref{eq7}).\fd 

\vskip0,3cm
Let $P_{\mathscr{A}}:\mathscr{A}(\R^n,E^*)\longrightarrow(\mathscr{A}(\R^n))^k$  be the continuous linear map defined by 
$$P_{\mathscr{A}}(f)=(\langle f,P_1\rangle,\ldots,\langle f,P_k\rangle),$$
and ${}^t\hskip-0,05cm P_{\mathscr{A}}$ its transpose. Let $(e_i)_{1\leq i\leq d}$ be a basis of $E$ and $(e_i^\ast)_{1\leq i\leq d}$ its dual basis. For $1\leq l\leq k$, we have  $P_l=\sum_{i=1}^dP_{li}\,e_i$ where $(P_{li})_{1\leq i\leq d}$ are homogeneous polynomials  with the same degree of homogeneity $s_l$ as $P_l$. 

For all function $f$ in $\mathscr{A}(\R^n,E^*)$, there exist $f_1,\ldots,f_d\in\mathscr{A}(\R^n)$ such that  $f=\sum_{i=1}^df_{i}\,e_i^\ast$. 
Then the spaces $\mathscr{A}(\R^n,E^*)$ and $(\mathscr{A}(\R^n))^d$ can be identified by the map $\sum_{i=1}^df_{i}\,e_i^\ast\longmapsto(f_1,\ldots,f_d)$. 

Let $P_{\mathscr{A}}$ denotes the map 
$$(f_1,\ldots,f_d) \longmapsto  (\sum_{i=1}^df_i\;P_{1i},\ldots,\sum_{i=1}^df_i\;P_{ki})
    ,$$
from $(\mathscr{A}(\R^n))^d$ to $(\mathscr{A}(\R^n))^k$.

Then we have the following lemma.
\begin{lem}\label{lem5}
The image of $\;{}^tP_{\mathscr{A}}$ is equal to ({\rm Ker}$P_{\mathscr{A}})^\circ$, the orthogonal space of {\rm Ker}$P_{\mathscr{A}}$.
\end{lem}
\pre\\
The inclusion $\im({}^tP_{\mathscr{A}})\subset({\rm Ker}P_{\mathscr{A}})^\circ$ is clear. When $\mathscr{A}(\R^n,E^*)$ is a Fréchet space it suffices to prove that the image of $P_{\mathscr{A}}$ is closed. Indeed, this implies that the image of $P_{\mathscr{A}}$ is a Fréchet space. Then Theorem 37.2 of \cite{TRE} insures that $\im({}^tP_{\mathscr{A}})$ is closed in $\mathscr{A}'(\R^n,E)$. 
On the other hand, Proposition 35.4 of \cite{TRE} implies that the closure of $\im({}^tP_{\mathscr{A}})$ is equal to $({\rm Ker}P_{\mathscr{A}})^\circ$.

Now we will prove that the image of $P_{\mathscr{A}}$ is closed. 

The case $\mathscr{A}(\R^n,E^*)=\mathcal{C}^\infty(\R^n,E^*)$ is a consequence of Theorem 0.1.1 of \cite{BSC}.

For $\mathscr{A}(\R^n,E^*)=\mathscr{S}(\R^n,E^*)$, we consider the map $\widehat{P}:(\mathcal{C}_{\overline{B}(0,1)}^\infty(\R^n))^d\longrightarrow(\mathcal{C}_{\overline{B}(0,1)}^\infty(\R^n))^k$ defined by  $$\widehat{P}=I_k\circ P_{\mathscr{S}}\circ J_d\;,$$ where $I_k$ is the topological isomorphism  from $(\mathscr{S}(\R^n))^k$ to $(\mathcal{C}_{\overline{B}(0,1)}^\infty(\R^n))^k$ defined by :
$$I_k:(f_1,\ldots,f_k)\longmapsto (I(f_1),\ldots,I(f_k)),$$ 
and  $J_d$ the topological isomorphism from $(\mathcal{C}_{\overline{B}(0,1)}^\infty(\R^n))^d$ to $(\mathscr{S}(\R^n))^d$ defined by :$$J_d:(f_1,\ldots,f_d)\longmapsto (I^{-1}(f_1),\ldots,I^{-1}(f_d)).$$ 

Then it is sufficient  to prove that the image of $\widehat{P}$ is closed. For  $f=(f_1,\ldots,f_d)\in(\mathcal{C}_{\overline{B}(0,1)}^\infty(\R^n))^d$, we have  
 $$\widehat{P}((f_1,\ldots,f_d))=(\widehat{P}_1(f),\ldots,\widehat{P}_k(f)),$$
   where, for all $1\leq j\leq k$, $$\widehat{P}_j(f)(x)=\left\{\begin{array}{ll}
 \displaystyle\sum_{i=1}^d f_i(x)\;\frac{P_{ji}(x)}{(1-||x||^2)^{\frac{s_j}{2}}}&\textrm{ if }x \in B(0,1)\\
    0 & \textrm{ otherwise }
    \end{array}.\right.$$
Let  $B:(\mathcal{C}_{\overline{B}(0,1)}^\infty(\R^n))^d\longrightarrow(\mathcal{C}_{\overline{B}(0,1)}^\infty(\R^n))^k$ defined by 
  $$B((g_1,\ldots,g_d))=(B_1(g),\ldots,B_k(g)),$$ where $$B_j(g)(x)=\sum_{i=1}^d g_i(x)\;P_{ji}(x),\quad\textrm{ for all } 1\leq j\leq k.$$
 Let 
$A:(\mathcal{C}_{\overline{B}(0,1)}^\infty(\R^n))^k\longrightarrow(\mathcal{C}_{\overline{B}(0,1)}^\infty(\R^n))^k$ defined by 
$$A((g_1,\ldots,g_k))=(A_1(g),\ldots,A_k(g)),$$ where
$$A_j(g)(x)=\left\{\begin{array}{ll}
 \displaystyle \frac{g_j(x)}{(1-||x||^2)^{\frac{s_j}{2}}}&\textrm{ if }x \in B(0,1)\\
    0 & \textrm{ otherwise }
    \end{array},1\leq j\leq k.\right.$$
It is clear that $A$ is a topological isomorphism, and that $\widehat{P}=A\circ B$. Again, by Theorem 0.1.1 of \cite{BSC}, the image of the map $B$ is closed. Consequently, the image of $\widehat{P}$ is closed. 

Now we will prove the result for $\mathscr{A}(\R^n,E^*)=\mathcal{C}^\infty_c(\R^n,E^*)$ which is not a Fréchet space. Let $m\in\N^\ast$. We consider the map 
$$P_m:f\longmapsto(\langle f,P_1\rangle,\ldots,\langle f,P_k\rangle),$$
from $\mathcal{C}_{\overline{B}(0,m)}^\infty(\R^n,E^*)$ to ($\mathcal{C}_{\overline{B}(0,m)}^\infty(\R^n))^k$.
We will  prove that $({\rm Ker}P_{\mathscr{A}})^\circ\subset\im({}^t\hskip-0,02cm P_{\mathscr{A}})$. Let $T$ be in $({\rm Ker}P_{\mathscr{A}})^\circ$. It is clear that the restriction of $T$  to the space $\mathcal{C}^\infty_{\overline B(0,m)}(\R^n,E^*)$ is in $({\rm Ker}P_m)^\circ$. Since  $\mathcal{C}^\infty_{\overline B(0,m)}(\R^n,E^*)$ is a Fréchet space, it follows, as in the first part of the proof, that $\im({}^t\hskip-0,02cm P_m)=({\rm Ker}P_m)^\circ$. So there exists a continuous linear form $\theta_m$ on $(\mathcal{C}^\infty_{\overline B(0,m)}(\R^n))^k$ such that $T$ is equal to ${}^t\hskip-0,02cm P_m(\theta_m)$  on $\mathcal{C}^\infty_{\overline B(0,m)}(\R^n,E^*)$. By the Hahn-Banach theorem we extend $\theta_m$ to a continuous linear form on $(\mathcal{C}^\infty_c(\R^n))^k$ denoted also by $\theta_m$. For all $m\in\N^\ast$ let $$C_m=\{x\in\R^n,\;m-2<||x||<m-\frac{1}{2}\}.$$ We fix a partition of unity $(\chi_m)_{m\in\N^\ast}$ subordinated to the covering $(C_m)_{m\in\N^\ast}$ of $\R^n$ with $\chi_m\in\mathcal{C}^\infty_c(\R^n)$, and we consider the distribution $\theta=\sum_{m=1}^\infty\chi_m\;\theta_m$. 

Let $f\in\mathcal{C}^\infty_c(\R^n,E^*)$.  We have 
$$\langle\theta\circ P_{\mathscr{A}},f\rangle=\sum_{m\in\N^\ast}\langle \theta_m,\chi_m\;P_{\mathscr{A}}(f)\rangle=\sum_{m\in\N^\ast}\langle \theta_m,P_{\mathscr{A}}(\chi_m\;f)\rangle.$$ Since $\supp\;\chi_m\subset\overline{B}(0,m)$, then 
$$\langle\theta\circ P_{\mathscr{A}},f\rangle=\sum_{m\in\N^\ast}\langle \theta_m,P_m(\chi_m\;f)\rangle\; \textrm{ and } \;\langle \theta_m,P_m(\chi_m\;f)\rangle=\langle T,\chi_m\;f\rangle.$$ Consequently, $$\langle\theta\circ P_{\mathscr{A}},f\rangle=\langle T,\sum_{m\in\N^\ast}
\chi_m\;f\rangle=\langle T,f\rangle.$$ So $T={}^t\hskip-0,02cm P_{\mathscr{A}}(\theta)$. 
\fd
 \begin{lem}\label{lem6}
The space  $\mathscr{A}'(\R^n,E)^G$ is contained in  ${}^t\hskip-0,03cm P_{\mathscr{A}}((\mathscr{A}(\R^n)^k)')$.
       \end{lem}
\pre\\
 By Lemma \ref{lem4} we have $\mathscr{A}'(\R^n,E)^G=\overline{\mathcal{C}^{\infty}_c(\R^n,E)^G}$     and
 by Lemma \ref{lem5}, $({\rm Ker}P_\mathscr{A})^\circ={}^t\hskip-0,03cm P_{\mathscr{A}}((\mathscr{A}(\R^n)^k)')$.   To conclude, we need to prove the inclusion  $\mathcal{C}^{\infty}_c(\R^n,E)^G\subset({\rm Ker}P_\mathscr{A})^\circ$. 

Let $\varphi\in\mathcal{C}^{\infty}_c(\R^n,E)^G$ and $f\in {\rm Ker} P_\mathscr{A}$. Then $\langle f(x),P_i(x)\rangle=0$ for all $x\in\R^n$ and all $1\leq i\leq k$.
By Theorem \ref{th2}, there exists a family of functions $\varphi_1,\ldots,\varphi_k\in\mathcal{C}^{\infty}_c(\R^n)^G$ such that  $$\varphi=\sum_{i=1}^k\varphi_i\;P_i.$$
 Then, $$\langle T_\varphi,f\rangle=\int_{\R^n}\langle f(x),\sum_{i=1}^k{\varphi_i(x)\;P_i(x)}\rangle\;dx=0.$$
\fd
\noindent
\preu{\bf\ref{th3}}\\
Let $\theta$ be in $\mathscr{A}'(\R^n,E)^G$. By Lemma \ref{lem6}, there exists $T\in(\mathscr{A}(\R^n)^k)'$ such that ${}^t\hskip-0,02cm P_{\mathscr{A}}(T)=\theta$. Then there is a family of distributions $(T_i)_{1\leq i\leq k}$ in $\mathscr{A}'(\R^n)$ such that  $$\langle T,(f_1,\ldots,f_k)\rangle=\sum_{i=1}^k \langle T_i,f_i\rangle,\quad\textrm{ for all }f_1,\ldots,f_k\in\mathscr{A}(\R^n),$$ and $\theta=\sum_{i=1}^kT_i\;P_i$. 

For $1\leq i\leq k$, we define, as for the $G$-covariant distributions in the proof of Lemma \ref{lem4}, the distribution $\theta_i$ in $\mathscr{A}'(\R^n)^G$  by 
\begin{equation}\label{eq8}\langle\theta_i,\varphi\rangle=\int_G{\langle T_i,\varphi\circ g\rangle\;d\mu(g)},\quad\textrm{ for all }\varphi\in\mathscr{A}(\R^n).\end{equation} 
Let $f$ be in $\mathscr{A}(\R^n,E^\ast)$. Since $\theta$ is $G$-covariant, then $$\langle\theta,f\rangle=\int_G\langle\theta,g\hskip-0,005cm\cdot\hskip-0,007cm  f\rangle\;d\mu(g).$$ So we have  $$\langle\theta,f\rangle=\sum_{i=1}^k{\int_G{\langle T_i,\langle g\hskip-0,005cm\cdot\hskip-0,007cm  f,P_i\rangle\rangle\;d\mu(g)}}.$$
The polynomials $(P_i)_{1\leq i\leq k}$ are $G$-covariant, then
$$\langle T_i,\langle g\hskip-0,005cm\cdot\hskip-0,007cm  f,P_i\rangle=\langle T_i,\langle f\circ g, P_i\circ g\rangle\rangle.$$
This implies by  Formula (\ref{eq8}) that  $$\langle\theta,f\rangle=\sum_{i=1}^k{\langle\theta_i,\langle f,P_i\rangle\rangle}=\langle\sum_{i=1}^k\theta_i\;P_i,f\rangle.$$ So $\theta=\sum_{i=1}^k\theta_i\;P_i$.\fd

\end{document}